\documentclass[11pt]{article}
\usepackage{latexsym}
\usepackage{amsfonts}
\makeatletter

\def\@footnotetext#1{\insert\footins{%

\footnotesize

 \interlinepenalty\interfootnotelinepenalty

 \splittopskip\footnotesep

 \splitmaxdepth \dp\strutbox \floatingpenalty \@MM

 \hsize\columnwidth \@parboxrestore

 \edef\@currentlabel{\csname p@footnote\endcsname\@thefnmark}\@makefntext
 {\rule{\z@}{\footnotesep}\ignorespaces
#1\strut}}}

\def\abstract{\small\quotation{\hskip-\parindent\sc Abstract.}}
\def\classification{\@ifnextchar [{\@xfootnotenext}%
 {\begingroup\let\protect\noexpand
 \xdef\@thefnmark{}\endgroup
 \@footnotetext}}

\title {}

\begin{document}

\classification {{\it 2000 Mathematics Subject Classification:} Primary 14E09, 14E25; Secondary 14A10, 13B25.\\
$\ast$) Partially supported by RGC Grant Project HKU 7134/00P.}

\begin{center}

{\bf \Large  Affine varieties  with equivalent  cylinders}

\bigskip

{\bf Vladimir Shpilrain}

\medskip

 and

\medskip

 {\bf Jie-Tai Yu}$^{\ast}$

\end{center}

\medskip

\begin{abstract}
\noindent    A well-known {\it cancellation problem} asks when, for two 
algebraic varieties $V_1, V_2 \subseteq {\bf C}^n$, the isomorphism 
of the cylinders $V_1 \times {\bf C}$ and  $V_2 \times {\bf C}$ implies 
the isomorphism of $V_1$ and  $V_2$.  

In this paper, we address   a related  problem: when 
 the {\it equivalence} (under an automorphism of ${\bf C}^{n+1}$) 
of two cylinders $V_1 \times {\bf C}$ and  $V_2 \times {\bf C}$ 
implies the equivalence of their bases $V_1$ and $V_2$ 
under an automorphism of  ${\bf C}^n$? 
We concentrate here on hypersurfaces  
and show that this problem establishes a strong connection between the 
Cancellation conjecture  of Zariski and the  Embedding conjecture of 
Abhyankar  and  Sathaye. We settle the problem   for 
a large class of polynomials. On the other hand, we give examples 
of equivalent cylinders with inequivalent bases (those cylinders, however, 
are not  hypersurfaces). 

 Another result of interest is that, for an arbitrary field $K$, 
 the equivalence of two polynomials in $m$
 variables under an automorphism of $K[x_1,..., x_n], ~n \ge m,$
~implies their equivalence under a {\it   tame}  automorphism 
of $K[x_1,..., x_{2n}]$. 

\end{abstract} 

\bigskip

\noindent {\bf 1. Introduction }

\bigskip

 Let ${\bf C}[x_1,..., x_n]$ be the polynomial algebra in $n$ variables
over the field ${\bf C}$. Any collection of polynomials $p_1,...,p_m$ from
${\bf C}[x_1,..., x_n]$ determines an algebraic variety
${\{}p_i=0, ~i=1,...,m{\}}$ in the affine space ${\bf C}^n$.
We shall denote
this algebraic variety by $V(p_1,...,p_m)$.
 \smallskip

 We say that two algebraic varieties $V(p_1,...,p_m)$ and $V(q_1,...,q_k)$
 are {\it isomorphic} if the algebras of residue
classes ${\mathbf C}[x_1,..., x_n]/\langle p_1,...,p_m \rangle$ and
${\mathbf C}[x_1,..., x_n]/\langle q_1,...,q_k \rangle$ are
isomorphic.
Here $\langle p_1,...,p_m \rangle$ denotes the ideal of
 ${\mathbf C}[x_1,..., x_n]$ generated by $p_1,...,p_m$.
 Thus, isomorphism that we consider here is algebraic, not geometric, 
 i.e., we actually consider isomorphism of what is called  {\it affine schemes}. 
 
 On the other hand, we say that two algebraic varieties $V(p_1,...,p_m)$
 and $V(q_1,...,q_k)$ are
{\it equivalent} if there is an automorphism of ${\bf C}^n$ that takes one
of them onto the other. Algebraically, this means there is 
an automorphism of ${\bf C}[x_1,..., x_n]$ that takes the ideal 
$\langle p_1,...,p_m \rangle$ to the ideal $\langle q_1,...,q_k \rangle$.

 Furthermore, a variety equivalent to  $V \times {\bf C}$ is called a 
{\it cylinder};
 a variety of the form $\{p=0\}$ is called a  {\it hypersurface}, and a 
hypersurface equivalent to $\{x_1=0\}$ is called a  {\it hyperplane}. 
In particular, a {\it cylindrical hypersurface} in ${\bf C}^n$ is  a variety 
of the form $\{p(x_1,...,x_m)=0\}$, where $m < n$.  
 \smallskip

 Below we list 3 conjectures relevant to the subject of the present paper. 
   The first two of them are very well known. 

\medskip 

\noindent {\bf Conjecture 1.}  (The  Cancellation conjecture  of Zariski). 
   Let $V \subseteq {\bf C}^n$ be a hypersurface. If $V \times {\bf C}$ is isomorphic to 
a hyperplane in ${\bf C}^{n+1}$, then $V$ is isomorphic to 
a hyperplane in ${\bf C}^n$. 
\medskip

\noindent  {\bf Conjecture 2.}  (The Embedding conjecture of 
Abhyankar  and  Sathaye).  If a hypersurface $V(p)$ in 
${\bf C}^n$ is isomorphic to the  
hyperplane  $V(x_1)$, then it is  equivalent to it.  

\medskip 

\noindent {\bf Conjecture 3.}    Let $V \subseteq {\bf C}^n$ be a hypersurface. 
If $V \times {\bf C}$ is equivalent to 
a hyperplane in ${\bf C}^{n+1}$, then $V$ is equivalent to 
a hyperplane in ${\bf C}^n$. Or, in purely algebraic language: 
if $p=p(x_1,..., x_n)$ ~and   $~\varphi(p) = x_1~$ for some automorphism
 $~\varphi$ of ${\bf C}[x_1,..., x_{n+1}]$, then  also $~\alpha(p) = x_1~$ 
for some automorphism $~\alpha$ of ${\bf C}[x_1,..., x_n]$. 
\medskip

   These conjectures are all unsettled in general, but there are some partial 
results. Conjecture 2 is settled for $n=2$ in \cite{AM}; 
 Conjecture 1 is 
settled for $n=2$ in \cite{AEH} and \cite{Miy1}, and for $n=3$ in \cite{Fujita}.
 In  \cite{Fujita}, Fujita, based on the work of 
Miyanishi  and Sugie  \cite{MiySugie}, proved that $V \times {\bf C}^n$ 
isomorphic to ${\bf C}^{n+2}$ implies that $V$ is 
isomorphic to ${\bf C}^2$. In \cite{Asanuma}, Asanuma gives an  idea for  
constructing a (possible) counterexample to Conjecture 1 over the field of 
reals.

For more information, see the survey \cite{Kraft}. 
\smallskip

  The focus of the present paper is on Conjecture 3. It   looks   similar 
to Conjecture 1; the   difference is that isomorphism here 
is replaced by equivalence. 

 All the three conjectures seem to be related; however, it is not so easy
to pinpoint precise relations between them. Here we prove the following 

\medskip

\noindent {\bf Proposition 1.1.}   

\noindent {\bf (a)} Conjectures 1 and  2 together imply Conjecture 3.
\smallskip

\noindent {\bf (b)} Conjectures 3 and  2 together imply Conjecture 1. 
\smallskip

\noindent {\bf (c)}  Conjecture 3 is true for  $n=2$ and  $3$. 

\medskip

 Thus, Conjecture 2 is, in a way, the  strongest of the three; in fact, 
since   Conjecture 3 is most likely true for any $n$, it is probably 
the case that Conjecture 2 just implies Conjecture 1. 

 Our proof of part  (c) in the case $n=3$ is based on deep results cited above, 
whereas the case $n=2$ can also be handled by somewhat more elementary methods, 
based on  a result of  \cite{CYu}. 
\smallskip

 In the course of our work on Conjecture 3, we were able to prove the following
result of independent  interest. Before we give 
the statement, we need a few more 
\medskip

\noindent {\bf   Definitions.} Let $K$ be  an arbitrary  
field. 
 An automorphism of $K[x_1,..., x_n]$ 
is called  {\it elementary} if it fixes all variables but one, say, $x_i$, 
 and maps $x_i$ to $x_i + f(x_1,...,x_{i-1}, x_{i+1},..., x_n)$. 
An automorphism of $K[x_1,..., x_n]$ 
is called  {\it tame} if it is a product of elementary  and linear 
automorphisms. A polynomial $p$ 
is called a {\it coordinate polynomial} (or simply a {\it coordinate}) if 
there is an 
  automorphism of $K[x_1,..., x_n]$ that takes $p$ to $x_1$.
A polynomial $p \in K[x_1,..., x_n]$ is called a {\it tame coordinate} 
 if it is a component of a tame  automorphism, i.e., if there is a 
tame  automorphism of $K[x_1,..., x_n]$ that takes $p$ to $x_1$. 
Finally, a polynomial $p \in K[x_1,..., x_n]$ is called a
 {\it stably tame coordinate} if it is a component of a tame  automorphism 
of $K[x_1,..., x_N]$  for some $N \ge n$. More generally, two 
polynomials $p, q \in K[x_1,..., x_n]$ are called 
 {\it stably (tame) equivalent} if $~\alpha(p) = q$ for some (tame) 
automorphism $~\alpha$ of $K[x_1,..., x_N]$, ~$N \ge n$.
\medskip

 We were able to prove: 
 \medskip

\noindent {\bf Theorem 1.2.} If two polynomials 
$p, q \in K[x_1,..., x_n]$ are equivalent, then they are 
stably tame equivalent. 
 \medskip

 In particular:

\medskip

\noindent {\bf Corollary 1.3.} Every coordinate is stably tame. 
\medskip

 Thus, to settle Conjecture 3, it is sufficient to consider cylinders 
equivalent under a tame automorphism. 
\smallskip

  Conjecture 3 also motivates the following generalization: 
\medskip

\noindent {\bf Problem 1.} Are any two stably equivalent polynomials 
equivalent? Or, in geometric language:  is it true that equivalence 
of two cylindrical  hypersurfaces $V(p) \times {\bf C}$ and  
$V(q) \times {\bf C}$ in ${\bf C}^{n+1}$ implies 
equivalence of their bases $V(p)$ and  $V(q)$ in  ${\bf C}^n$ ? 

\medskip 

 A weaker version of this problem is still quite  hard: 
\medskip

\noindent {\bf Problem 2.} Is it true that equivalence 
of two cylindrical  hypersurfaces $V(p) \times {\bf C}$ and  
$V(q) \times {\bf C}$ in ${\bf C}^{n+1}$ implies 
isomorphism  of their bases $V(p)$ and  $V(q)$ in  ${\bf C}^n$ ? 
\medskip

 A generalization of Conjecture 3 in another direction is 

\medskip

\noindent {\bf Problem 3.}  
 Can a cylinder in ${\bf C}^n$ 
be  equivalent to a non-cylinder  in ${\bf C}^n$ under an automorphism 
of ${\bf C}^{m}, ~m > n$ ? 

\medskip 

Problem 1 is, in our opinion, the most interesting one, and it probably 
holds the key to many mysteries of affine algebraic geometry. 
  
  We  were able to prove that the  answer to 
 Problem 1 is affirmative for two rather large classes of polynomials 
(Proposition  1.4 and Theorem  1.5 below).   The first result is just 
a simple observation: 
\medskip 

\noindent {\bf Proposition  1.4.} Let $p=p(x_1,...,x_n)$ be 
of the form $x_1^{M_1}\cdot ...\cdot x_n^{M_n} + 
\sum_j c(j) \cdot x_1^{i_1(j)}...x_n^{i_n(j)}$, 
~where $c(j)$ are some coefficients,  all $M_i >0$, all monomials under 
the sum are different from $x_1^{M_1}\cdot ...\cdot x_n^{M_n}$, 
 and, for every $k, j$, one has $i_k(j) \le M_k$. 
  Then  any  polynomial $q=q(x_1,...,x_n)$ stably equivalent to $p$, 
is equivalent to it. 
 Or, in geometric language:  for any  hypersurface $V(q)$ in ${\bf C}^n$,   
 equivalence of \ $V(p) \times {\bf C}$ and  
$V(q) \times {\bf C^k}$ in ${\bf C}^{n+k}$  implies 
equivalence of  $V(p)$ and  $V(q)$ in  ${\bf C}^n$. 
\medskip

  We note that if $p$  is a coordinate in 
${\bf C}[x_1,..., x_n]$, then  
$p$ cannot possibly satisfy the condition of Proposition  1.4 since, 
by a result of Hadas \cite{Hadas},  
all vertices of the Newton polytope of a coordinate polynomial must be 
on coordinate hyperplanes. 
 Therefore,  Proposition  1.4  does not settle Conjecture 3.

 On the other hand, it should be pointed out that, even if a given 
polynomial $p$ does not satisfy the condition of Proposition  1.4, there 
might be an automorphic image of $p$ that does, and then the result will 
hold for $p$ as well. For example, $p(x,y)=(x+y)^2 y^2$ does not have the form 
required in Proposition  1.4, but  it has an automorphic image  
 $q(x,y)=x^2 y^2$ that does. 
We discuss the two-variable case in more detail 
in Section 2, after the proof of Proposition  1.4. 
\smallskip

 Another,  more interesting, class of polynomials for which 
Problem 1 has the affirmative answer, is given by the following 
\medskip

\noindent {\bf Theorem  1.5.} Suppose a polynomial $p=p(x_1,...,x_n)$ has 
the following property: if $\psi(p)=p$ for some injective polynomial mapping 
$\psi$  of ${\bf C}[x_1,..., x_n]$, then $\psi$ must be an automorphism 
of ${\bf C}[x_1,..., x_n]$. Then any  polynomial $q=q(x_1,...,x_n)$ stably 
equivalent to $p$, is equivalent to it. 
\medskip 

 Polynomials that satisfy the condition of Theorem  1.5 are called 
{\it test polynomials} for injective polynomial mappings -- see \cite{ESh}. 
Jelonek \cite{Jelonek} showed 
that a generic polynomial of degree
$> n$   in ${\bf C}[x_1,..., x_n]$ is a test 
polynomial for injective mappings. (The statement ``a generic polynomial  
of degree $d$ has a property W" means that there exists
a Zariski open dense subset $U\subseteq {\cal H}_d$ of the set 
 ${\cal H}_d$  of all polynomials of degree $d$,
such that every element of $U$ has the property W.) Thus, 
 the class of polynomials covered by Theorem 1.5 is really 
large. Still, it does not include any coordinate polynomials, as was 
observed in \cite{ESh}. 

\medskip 

 Then, we were able to answer Problem 3 in the affirmative: 

\medskip

\noindent {\bf Theorem  1.6.} For any  $k \ge 1$ and  $n \ge k+2$, there is 
a variety $V \times {\bf C}^k$ in ${\bf C}^n$ and a  non-cylinder $U$
  in ${\bf C}^n$, such that $U$ is equivalent to 
$V \times {\bf C}^k$ in ${\bf C}^{2n}$.

\medskip

 The simplest example (that we know) illustrating  Theorem  1.6 would be 
$U= \\
\{x(1+xy+z^2)=1\}$ in ${\bf C}^3$ and  $V =\{xy=1\}$ in ${\bf C}^2.$ 
The varieties 
$U$ and $V \times {\bf C}$ are inequivalent in ${\bf C}^3$, but 
are equivalent in ${\bf C}^6$ (in fact, they are equivalent even  in ${\bf C}^4$). 
\medskip 

 We also mention here an example (due to Danielewski \cite{Dan}, 
unpublished) of {\it isomorphic}  cylindrical hypersurfaces  
$V(p) \times {\bf C}$ and  $V(q) \times {\bf C}$ with non-isomorphic bases 
$V(p)$ and  $V(q)$ in ${\bf C}^3$. In his example, $p=p(x,y,z)=xy-z^2+1; 
~q=q(x,y,z)=x^2y-z^2+1$.

Recently,  Bandman and   Makar-Limanov  \cite{BML} uncovered a geometric
 reason why two cylinders with non-isomorphic bases can possibly be isomorphic. 
 In our concluding Section 3, we give an explicit  algebraic isomorphism 
(due to P. Russell) for the cylinders in Danielewski's example,   to show 
how complicated it is. This also establishes the fact that 
Danielewski's cylinders are isomorphic over any ground 
field   of characteristic $0$, not just over ${\bf C}$. 

  Furthermore,  a 
combination of Danielewski's example with our characterization of isomorphic 
varieties \cite[Corollary 1.2]{ShYu2} yields an example of {\it equivalent} 
cylinders with inequivalent (even non-isomorphic !) 
bases; those cylinders, however, are not  hypersurfaces: 
\medskip

\noindent {\bf Proposition 1.7.}  As in Danielewski's example, let 
$p(x,y,z)=xy-z^2+1; ~q(x,y,z)=x^2y-z^2+1$. Then the varieties 
$V_1=V(p(x,y,z), t, u, v, w)$ and $V_2=V(q(x,y,z), t, u, v, w)$ in ${\bf C}^7$ 
are not isomorphic, whereas the cylinders 
$V_1 \times {\bf C}$ and  $V_2 \times {\bf C}$ are equivalent 
in ${\bf C}^8$.  \\

\noindent {\bf 2. Proofs }

\bigskip

{\bf  Proof of Proposition 1.1.}
 We start with    part {\bf (a)}. We are  going to prove this 
statement in the following form: 
 if both Conjectures 1 and  2 are true for any $n$, 
then for any $k \ge 1$, whenever 
$V \times {\bf C}^k$ is equivalent to 
a hyperplane in ${\bf C}^{n+k}$, one has $V$  equivalent to 
a hyperplane in ${\bf C}^n$.

 Let $V=V(p), ~p=p(x_1,...,x_n)$, and suppose that $p$ is a coordinate 
in ${\bf C}[x_1,..., x_N]$, ~$N=n+k$. Then, in particular, the hypersurface 
$\{p=0\}$ in ${\bf C}^N$ is isomorphic to ${\bf C}^{N-1}$. If Conjecture 1  is true, 
this implies that $\{p=0\}$ in ${\bf C}^n$ is isomorphic to ${\bf C}^{n-1}$.
Now if  Conjecture 2 is true, this implies that $p$ is a coordinate 
in ${\bf C}[x_1,..., x_n]$, hence $V$ is equivalent to 
a hyperplane in ${\bf C}^n$. $\Box$ 

\medskip

 For part {\bf (b)}, let $V=V(p), ~p=p(x_1,...,x_n)$, and suppose that 
the hypersurface 
$\{p=0\}$ in ${\bf C}^N, ~N>n,$ is isomorphic to ${\bf C}^{N-1}$. 
If Conjecture  2 is true, this implies that $p$ is a coordinate 
in ${\bf C}[x_1,..., x_N]$. Then Conjecture 3 implies 
that $p$ is a coordinate in ${\bf C}[x_1,..., x_n]$ as well. 
In particular, the hypersurface 
$\{p=0\}$ in ${\bf C}^n$ is isomorphic to ${\bf C}^{n-1}$.

\medskip

  For part {\bf (c)}, we start with $n=2$. Let $V=V(p), 
~p=p(x_1,x_2)$, and suppose that for some $N \ge 3$, $p$ is a coordinate 
in ${\bf C}[x_1,..., x_N]$. Then, in particular, the hypersurface 
$\{p=0\}$ in ${\bf C}^N$ is isomorphic to ${\bf C}^{N-1}$.
Then, by the result of Abhyankar-Eakin-Heinzer 
\cite{AEH} and   Miyanishi \cite{Miy1} cited in the Introduction, 
the hypersurface $\{p=0\}$ in ${\bf C}^2$ is isomorphic to ${\bf C}$. 
This implies, by the Abhyankar-Moh theorem \cite{AM}, that $p$ is a coordinate 
in ${\bf C}[x_1,x_2]$. 
\smallskip

 Now we get to $n=3$. The proof here is similar. Let $V=V(p), 
~p=p(x_1,x_2, x_3)$. If $p$ is a coordinate 
in ${\bf C}[x_1,..., x_N]$ for some  
$N > 3$, then also for any $c \in {\bf C}$, ~$p-c$ is a coordinate 
in ${\bf C}[x_1,..., x_N]$. Therefore, for any $c \in {\bf C}$, 
the hypersurface 
$\{p=c\}$ in ${\bf C}^N$ is isomorphic to ${\bf C}^{N-1}$. 
 Then, again for any $c \in {\bf C}$, the hypersurface 
$\{p=c\}$ in ${\bf C}^3$ is isomorphic to ${\bf C^2}$ by
the result of Fujita \cite{Fujita}  cited in the Introduction.
 Now this implies, by a recent result of Kaliman \cite{K} 
(see also \cite{KZ}),  
that $p$ is a coordinate in ${\bf C}[x_1,x_2,x_3]$, 
hence $V$ is equivalent to 
a hyperplane in ${\bf C}^3$. $\Box$ 

\medskip

\noindent {\bf  Proof of Theorem 1.2.}  
 Let $\varphi$ be an automorphism   of $K[x_1,..., x_n]$. 
By just replacing all $x_i$ with $y_i$, we make a ``copy" of 
$\varphi$ acting on $K[y_1,..., y_n]$.   
We are going to show that the following automorphism  
 of $K[x_1,..., x_n, y_1,...,y_n]$ is tame:
 $\psi: x_i \to \varphi(x_i), 1 \le i \le n; ~y_i \to -\varphi^{-1}(y_i), 
1 \le i \le n$. 

 Let $\alpha: x_i \to x_i + \varphi(y_i), ~y_i \to y_i, ~1 \le i \le n$, 
~and  ~$\beta: x_i \to x_i, ~y_i \to y_i-\varphi^{-1}(x_i), ~1 \le i \le n$, 
be two automorphisms of $K[x_1,..., x_n, y_1,...,y_n]$. They both are 
obviously tame. 

 Now let 
$$\psi_1 = \alpha \circ \beta = \beta(\alpha): 
x_i \to x_i + \varphi(y_i-\varphi^{-1}(x_i))=\varphi(y_i), 
~y_i \to y_i-\varphi^{-1}(x_i).$$ 
Compose $\psi_1$ with a linear 
automorphism $\pi: x_i \to y_i; ~y_i \to x_i$. We get: 
$\psi_2 = \psi_1 \circ \pi = \pi(\psi_1): x_i \to \varphi(x_i), 
~y_i \to x_i-\varphi^{-1}(y_i)$. 

 Finally, compose $\psi_2$ with $\tau: x_i \to x_i, ~y_i \to y_i + \varphi(x_i)$ 
to get 
$$\psi = \tau(\psi_2): x_i \to \varphi(x_i), 
~y_i \to x_i-\varphi^{-1}(y_i + \varphi(x_i))= -\varphi^{-1}(y_i).$$

 This completes the proof. $\Box$ 

\medskip

\noindent {\bf Remark 2.1.} The same argument establishes a somewhat more 
general result: if $p_i, q_i \in K[x_1,..., x_n], ~1 \le i \le m$, 
and $\varphi(p_i)=q_i, ~1 \le i \le m$, for some automorphism $\varphi$ 
of $K[x_1,..., x_n]$, then $\psi(p_i)=q_i, ~1 \le i \le m$, for some 
tame automorphism $\psi$ of $K[x_1,..., x_{2n}]$. In particular, every 
{\it coordinate tuple} of polynomials in $K[x_1,..., x_n]$ is stably tame 
in the sense that it is part of a tame automorphism of $K[x_1,..., x_{2n}]$.
\medskip

\noindent {\bf Remark 2.2.} The same argument can be used, in fact, 
for many type of free algebras, thus establishing the same result, 
 in particular, for a free associative algebra   
$K\langle x_1,..., x_n \rangle$.  This also implies that 
every coordinate of $K[x_1,..., x_n]$ can be {\it lifted} 
to a primitive element (this is what coordinates are called in 
non-commutative setting) of $K\langle x_1,..., x_{2n}\rangle$.  

\medskip

\noindent {\bf  Proof of Proposition  1.4.} Let $\varphi(p)=q$ for some 
automorphism $\varphi$ of ${\bf C}[x_1,..., x_N], ~N>n$. Then 
all $x_i$ with $i>n$ have to cancel out in $\varphi(p)$. However, 
because of the presence of the ``dominating" monomial 
$x_1^{M_1}\cdot ...\cdot x_n^{M_n}$ in $p$, if there is some 
$x_i$ with $i>n$ in some of the $\varphi(x_j), ~1 \le j \le n$, then 
this $x_i$ will not cancel out in $\varphi(p)$. Therefore, 
none of the $\varphi(x_j), ~1 \le j \le n$, depends on any 
$x_i, ~i>n$, in which case the restriction of $\varphi$ to 
${\bf C}[x_1,..., x_n]$ must be an automorphism 
of ${\bf C}[x_1,..., x_n]$. $\Box$ 

\medskip

\noindent {\bf Remark 2.3.} As we have mentioned in the Introduction, 
even if a given 
polynomial $p$ does not satisfy the condition of Proposition  1.4, there 
might be an automorphic image of $p$ that does, and then the result will 
hold for $p$ as well. Here we note that, in the case where $p$ is 
a two-variable polynomial over ${\bf C}$, it is algorithmically possible 
to find out whether some automorphic image of $p$ has the form 
required by Proposition  1.4. Indeed, Wightwick   \cite{Wightwick} 
showed that, if some automorphic image of $p$ has lower degree than 
$p$ does, then there is a single elementary automorphism that reduces 
the degree of $p$. At the same time, it is easy to see that no 
elementary automorphism can reduce  the degree of a polynomial 
that satisfies the condition of Proposition  1.4. 
 Thus, given a polynomial $p$, we keep applying 
elementary automorphisms that reduce the degree, until the degree 
becomes irreducible. Then we check if there is a linear or 
elementary automorphism preserving the degree of the obtained polynomial, 
that produces a polynomial in the  form required by Proposition  1.4. 
 For more details, we refer to \cite{Wightwick}. 

\medskip

\noindent {\bf  Proof of Theorem  1.5.} Let $\varphi(p)=q$ for some 
automorphism $\varphi$ of ${\bf C}[x_1,..., x_N], ~N>n$. It will be 
sufficient to prove that there is an injective mapping $\psi_\varphi$ of 
${\bf C}[x_1,..., x_n]$ such that $\psi_\varphi(p)=q$. 
 Indeed, if this is the case, then, similarly, there is an injective 
mapping $\psi_{\varphi^{-1}}$ of 
${\bf C}[x_1,..., x_n]$ such that $\psi_{\varphi^{-1}}(q)=p$. 
 Then $\psi_{\varphi^{-1}}\psi_\varphi(p)=p$, and $\psi_{\varphi^{-1}}\psi_\varphi$ 
is injective. The result follows. 

 We also may assume, by our  Theorem 1.2, that the automorphism  $~\varphi$ 
is tame. 
 Thus, the proof of the theorem will be complete if we establish 
the following 
\smallskip

\noindent {\bf Proposition 2.4.}  Let $~\varphi$ be a tame automorphism of 
${\bf C}[x_1,..., x_{n+k}]$, ~$n, k \ge 1$,  and let 
$~\varphi_{q_1,...,q_k}$ be the 
restriction  to ${\bf C}[x_{i_1},..., x_{i_n}]$ of the  homomorphism  
obtained by 
replacing some $k$ variables $x_{j_1},..., x_{j_k}$ 
 with     polynomials  
$q_j$, that depend on other variables $x_{i_1},..., x_{i_n}$ only, 
 in every 
$\varphi(x_{i_j}), ~1 \le j \le n$. Then, for some 
 polynomials $q_j$,  the mapping 
$\varphi_{q_1,...,q_k}$ is injective on ${\bf C}[x_{i_1},..., x_{i_n}]$. 
\smallskip 

 We are going to give a proof here for $k=1$ since in the general case, 
the proof goes along exactly the same lines, but the notation becomes 
intractable.  Also for notational convenience, we will assume that 
$\varphi_q$ is the 
restriction  to ${\bf C}[x_1,..., x_{n}]$ of the  homomorphism  $~\varphi$. 
Thus, we are going to   give a proof of the following  
\smallskip

\noindent {\bf Proposition 2.4.1.}   Let $~\varphi$ be a tame automorphism of 
${\bf C}[x_1,..., x_{n+1}]$, ~$n \ge 1$,  and let $~\varphi_q$ be the 
restriction of $~\varphi$ to ${\bf C}[x_1,..., x_{n}]$ obtained by 
replacing $x_{n+1}$ with   a  polynomial $q=q(x_1,..., x_{n})$ in every 
$\varphi(x_i), ~1 \le i \le n$. Then, for some 
 polynomial $q$, the mapping 
$\varphi_q$ is injective on ${\bf C}[x_1,..., x_{n}]$. 
\smallskip

 It will be convenient to single out   one (obvious) preliminary  
statement:  
\smallskip

\noindent {\bf Lemma 2.5.} 
In the notation of Proposition 2.4.1,  
 let the mapping $\varphi_q$ be injective on 
${\bf C}[x_1,..., x_{n}]$. Let $\varphi_q: x_i \to h_i, ~1 \le i \le n$. 
Then, for any   $m \ge 0$   and $j, ~1 \le i \le n$, 
 the mapping $\varphi_{q'}$ is injective on 
${\bf C}[x_1,..., x_{n}]$, too, where $q'=q \cdot h_j^m$. 
\smallskip

\noindent {\bf Proof of Lemma 2.5.} This statement says that if polynomials 
$h_1, h_2,...,h_n, q$ are algebraically independent, then 
polynomials $h_1, h_2,...,h_n, q \cdot h_j^m$ are 
algebraically independent, too. This is fairly obvious from considering 
the quotient fields  ${\bf C}(h_1,..., h_{n}, q)$ and 
${\bf C}(h_1,..., h_{n}, q \cdot h_j^m)$. $\Box$ 

\smallskip

 The reason why we need this flexibility in  choosing the polynomial $q$, 
will be clear below.

\smallskip

\noindent {\bf Proof of Proposition 2.4.1.}   We use induction on the number 
of elementary  and linear automorphisms 
in a decomposition of $~\varphi$. Thus, we assume that for some particular 
$\varphi$, we have found a polynomial $q$ such that $\varphi_q$ 
is injective on ${\bf C}[x_1,..., x_{n}]$, and now we want to find such a 
$q$ for a composition of $~\varphi$ with an elementary or linear 
automorphism. A linear automorphism, in its turn, 
is a product of  elementary  automorphisms, permutations of variables,  and 
multiplications (of some of the variables) by non-zero constants. 
 Composing $~\varphi$ with an   automorphism of one of
 the latter two kinds does not present any difficulty. 
 Thus, we may assume 
that we are composing $~\varphi$ with an elementary  automorphism, call it 
$\epsilon$.  Here we have to consider two cases: 
\smallskip

\noindent {\bf (i)} $\epsilon$ fixes all 
variables except some $x_i$ with $1 \le i \le n$. We may as well assume that 
$i=1$. Let $\varphi: x_i \to g_i, ~1 \le i \le n+1$, and let 
$\epsilon: x_1 \to x_1+ f(x_2,...,x_{n+1}); x_i \to x_i, ~2 \le i \le n+1$. 
We claim that the mapping $(\epsilon(\varphi))_{q'}$ is injective on 
${\bf C}[x_1,..., x_{n}]$ for some  $q'$  of the form $q \cdot g_j^m$. 

By way of contradiction, assume that  for some 
polynomial  $p$, one has 
$$p(g_1(x_1+f(x_2,...,x_{n},q'), ~x_2,...,x_n,q'), ...,
g_n(x_1+f(x_2,...,x_{n},q'), ~x_2,...,x_n,q'))=0.$$
 On the other hand, by the inductive assumption, we must have 
$$p(g_1(x_1, x_2,...,x_n,q'), ...,g_n(x_1, x_2,...,x_n,q')) \ne 0.$$ 
 The former equality is obtained from the latter inequality by applying 
the following mapping of the algebra ${\bf C}[x_1,..., x_{n}]$: 
$x_1 \to x_1+ f(x_2,...,x_{n},q'), ~x_i \to x_i, ~2 \le i \le n$. 
This  mapping is injective since $x_1$ cannot cancel out in 
$x_1+ f(x_2,...,x_{n},q')$ -- this is where we use the 
flexibility in  choosing the polynomial $q'$ provided by Lemma 2.5. 
Upon multiplying $q$ by  an appropriate $g_j^m$, we make sure that either 
$f(x_2,...,x_{n},q')$ has no $x_1$ whatsoever, or  it has $x_1$ in a 
monomial of degree at least 2. 

 Thus, we have  got a contradiction, which completes
 the proof in this case. 
\smallskip

\noindent {\bf (ii)} $\epsilon$ fixes all variables except $x_{n+1}$. 
Let  $\epsilon: x_{n+1} \to x_{n+1}+ f(x_1,...,x_{n}); x_i \to x_i, 
~1 \le i \le n$. In this case, from the inductive assumption, 
it is obvious that, if we take 
$q'=q-f$, then the mapping $(\epsilon(\varphi))_{q'}$ is going to be injective
 on ${\bf C}[x_1,..., x_{n}]$. $\Box$ 

\smallskip 

 We conclude this section with the 

\smallskip 

\noindent {\bf Proof of Theorem 1.6.} First we give  a proof 
for $k=1, ~n=3$ to simplify the notation; then we explain how this 
proof generalizes easily to arbitrary $k \ge 1$ and  $n \ge k+2$. 

 We are going to show that $U= \{x(1+xy+z^2)=1\}$ is isomorphic 
to  $W =\{xy=1\}$ in ${\bf C}^3.$ The latter variety obviously is 
of the form $V \times {\bf C}$, whereas the former is not. 

As in \cite{ShYu2}, it  will be technically more 
convenient  to write   algebras of residue classes as 
``algebras with relations", i.e., for example, instead of 
${\bf C}[x_1,...,x_n]/\langle p(x_1,...,x_n) \rangle$ we shall write 
$\langle x_1,...,x_n \mid p(x_1,...,x_n)=0\rangle$. 

 Now we  get the following chain of ``elementary" isomorphisms:
\smallskip

\noindent 
  $\langle x, y, z \mid x(1+xy+z^2)=1 \rangle = 
\langle x, y, z \mid x(1+xy+z^2)=1, ~xy(1+xy+z^2)=y \rangle = \\
\langle x, y, z \mid x +x^2y+xz^2=1, ~xy+x^2y^2+xyz^2=y \rangle ~\cong \\ 
\langle x, y, z, u \mid x +x^2y+xz^2=1, ~xy+x^2y^2+xyz^2=y, ~u=xy \rangle
~\cong\\
 \langle x, y, z, u \mid x +xu+xz^2=1, ~y=u+u^2+uz^2, ~u=xy \rangle ~\cong \\
\langle x, z, u \mid x +xu+xz^2=1, ~u=x(u+u^2+uz^2) \rangle ~\cong \\
\langle x, y, z \mid x +xy+xz^2=1, ~y=x(y+y^2+yz^2) \rangle =
\langle x, y, z \mid x +xy+xz^2=1 \rangle = \\
\langle x, y, z \mid x(1+y+z^2)=1 \rangle ~\cong 
\langle x, y, z \mid xy=1 \rangle.$

\smallskip

 Thus, $U= \{x(1+xy+z^2)=1\}$ is isomorphic to  $W =\{xy=1\}$ in ${\bf C}^3.$ 
Upon     replacing $U$ with $U_m = \{x(1+xy+z_1^2+...+z_m^2)=1\}$, where 
$z_1,..., z_m$ are variables, and using the same chain of ``elementary" 
isomorphisms, we get examples, for any $m \ge 1$, of   non-cylinders $U_m$ 
isomorphic to  $W =\{xy=1\}$ in ${\bf C}^{m+2}$. Then, upon 
replacing  $U_m$ with $U_{m, r} = \{x \cdot t_1 \cdot ...\cdot t_r \cdot 
(1+xy+z_1^2+...+z_m^2)=1\}$   and   $W$ with $W_r=\{xy\cdot t_1 \cdot ...\cdot t_r=1\}$, 
where $t_1,..., t_r$ are variables, and using the same chain of  
isomorphisms, we get isomorphism of $U_{m, r}$ to $W_r$ in ${\bf C}^{m+2+r}$
 for any $m \ge 1$, ~$r \ge 0$. 

 To show that, for instance, the polynomial $x(1+xy+z_1^2+...+z_m^2)$ is not 
equivalent to any polynomial in less than $(m+2)$ variables, it is clearly 
sufficient to show the same for the polynomial $q=1+xy+z_1^2+...+z_m^2$. 
The gradient of the latter polynomial is $(y, ~x, ~2z_1, ..., ~2z_m)$. 
Now the ``chain rule" for partial derivatives shows that, if we apply an 
automorphism $\phi$ to $q$, then the gradient of $\phi(q)$ must contain 
all the variables $x, y, z_1, ..., z_m$. Then the same is true for 
the polynomial $\phi(q)$ itself. 

 Finally, to get the equivalence 
 claimed in the statement of Theorem 1.6 out of the above isomorphism, 
we just recall (see e.g. \cite[Corollary 1.2]{ShYu2}) that 
if two algebraic varieties 
$~V(p_1,...,p_m)$ ~and ~$V(q_1,...,q_k)$ ~in  ${\bf C}^n$ are 
isomorphic, then the varieties ~$V(p_1,...,p_m, x_{n+1},...,x_{2n})$
~and  
$~V(q_1,...,q_k, x_{n+1},...,x_{2n})$ in  ${\bf C}^{2n}$ are 
equivalent under a (tame) automorphism of ${\bf C}^{2n}$. 
$\Box$  \\

\noindent {\bf 3. Around Danielewski's example}
\bigskip

  In Danielewski's example \cite{Dan} 
of two isomorphic cylindrical hypersurfaces  
$V(p) \times {\bf C}$ and  $V(q) \times {\bf C}$ with non-isomorphic bases 
$V(p)$ and  $V(q)$ in ${\bf C}^3$, we have $p=p(x,y,z)=xy-z^2+1; ~q=q(x,y,z)=
x^2y-z^2+1$. We start this section by giving an explicit  algebraic 
isomorphism 
(due to P.Russell) for cylinders in Danielewski's example. We are 
grateful to  P.Russell for kindly permitting us to use his observation here.

\smallskip

\noindent {\bf Remark 3.1.} (Russell, unpublished). Let $p=p(x_1,y_1,z_1)=
x_1y_1-z_1^2+1; ~q=q(x_2,y_2,z_2)=x_2^2y_2-z_2^2+1$, and  let  $K$ be an 
arbitrary field   of characteristic $0$. 
The following 
mapping $\varphi$ from $K[x_1, y_1, z_1, u]$ 
 to $K[x_2, y_2, z_2, v]$  induces an isomorphism between 
algebras of residue classes $K[x_1, y_1, z_1, u]/\langle p \rangle$ and 
$K[x_2, y_2, z_2, v]/\langle q \rangle$: 
\smallskip

\noindent $\varphi : x_1 \to x_2; ~y_1 \to x_2v^2 + 2z_2v +x_2y_2; 
~z_1 \to x_2v + z_2; \\
u \to x_2v^3 + 3z_2v^2 +3x_2y_2v +  y_2z_2$. 
\smallskip

 It is easy to check that $\varphi$ induces a homomorphism 
between $K[x_1, y_1, z_1, u]/\langle p \rangle$ and 
$K[x_2, y_2, z_2, v]/\langle q \rangle$. 

We are going to show 
that $\varphi$ is {\it onto}. We already have 
$x_2$ (mod $\langle q \rangle$) in the image; a straightforward 
computation shows that $\varphi(y_1 z_1)- x_2 \cdot \varphi(u)=
-2x_2^2y_2v + 2z_2^2v=2v(z_2^2-x_2^2y_2)=2v$ (mod $\langle q \rangle$). 
 Thus, we have $v$ (mod $\langle q \rangle$), and therefore also 
$z_2$ (mod $\langle q \rangle$) in the image. Now inspection of 
$\varphi(y_1)$ shows that we have 
$x_2 y_2$ (mod $\langle q \rangle$),   hence also 
$x_2^2 y_2^2$ (mod $\langle q \rangle$) in the image. From 
$\varphi(u)$ we now see  that we have 
$y_2 z_2$ (mod $\langle q \rangle$),   hence also
$y_2 z_2^2$ (mod $\langle q \rangle$) in the image. This finally 
gives $y_2 z_2^2 - x_2^2 y_2^2= y_2$ (mod $\langle q \rangle$) 
in the image. 

 Thus, $\varphi$ is {\it onto}. If $\varphi$ were not one-to-one, 
then the algebra $K[x_2, y_2, z_2, v]/\langle q \rangle$ would be 
isomorphic to an algebra $K[x_1, y_1, z_1, u]/J$, 
where $J$ is an ideal that properly contains $\langle p \rangle$. 
This isomorphism then implies that $J$ can be generated (as 
an ideal) by a single polynomial.   Since $p$ is irreducible, 
this generating polynomial has to be $p$. 

 Therefore, $\varphi$ is one-to-one, hence an isomorphism. 
\smallskip

 Now we get to 
\smallskip

\noindent {\bf Proof of Proposition 1.7.} First of all, the varieties 
$V_1=V(p(x,y,z), t, u, v, w)$ and $V_2=V(q(x,y,z), t, u, v, w)$ in ${\bf C}^7$ 
are not isomorphic since if they were, then $V(p(x,y,z))$ and $V(q(x,y,z))$ 
would be isomorphic in ${\bf C}^3$ which is known not to be the case. 

 On the other hand, we know that the algebras of residue
classes ${\bf C}[x, y, z, t]/\langle p \rangle$ and 
${\bf C}[x, y, z, t]/\langle q \rangle$ are isomorphic. This implies, by 
\cite[Corollary 1.2]{ShYu2}, that the varieties 
$V(p(x_1, y_1, z_1), u, x_2, y_2, z_2, v)$ and $V(q(x_1,y_1,z_1), u, x_2, y_2, z_2, v)$ 
are equivalent under an automorphism of ${\bf C}^8$. We note that the latter 
equivalence also follows from a result of  Asanuma \cite{Asanuma}. $\Box$
\medskip

 Finally, we make one more observation inspired by Danielewski's example: 
\medskip

\noindent {\bf Proposition 3.2.} Let $y, x_1,...,x_m, ~m \ge 1,$ be variables, and 
 $p=p(x_1,...,x_m)$   an arbitrary polynomial. Then  hypersurfaces 
$y \cdot p = 1$ and $y \cdot p^k = 1$ are isomorphic in ${\bf C}^{m+1}$ 
for any $k \ge 1$. 
\medskip

    Speaking somewhat informally, adding     a polynomial $q(z_1,...,z_k)$ 
in new variables $z_i$, to both polynomials in the statement of 
Proposition 3.2, 
can ``spoil" this simple but delicate isomorphism (in the course of the proof, 
we shall see why), but, apparently, sometimes  
 this  isomorphism    ``survives" 
in higher dimensions, hence the isomorphism of cylinders.  Understanding 
when exactly this happens can be the key to constructing a 
counterexample to (some of) the conjectures  mentioned in this paper.  
\smallskip

\noindent {\bf Proof of Proposition 3.2.}  
Here we  get the following chain of ``elementary" isomorphisms:
\smallskip

\noindent  $\langle y, x_1,...,  x_m \mid y \cdot p = 1 \rangle = 
 \langle y, x_1,...,  x_m \mid y \cdot p = 1, y^2 \cdot p = y \rangle ~\cong \\
 \langle y, x_1,...,  x_m, u \mid y \cdot p = 1, y^2 \cdot p = y, u=y^2 \rangle ~\cong \\
 \langle y, x_1,...,  x_m, u \mid u \cdot p^2 = 1, u \cdot p = y, u=u^2p^2 \rangle ~\cong \\
\langle x_1,...,  x_m, u \mid u \cdot p^2 = 1, u=u^2p^2 \rangle =
\langle x_1,...,  x_m, u \mid u \cdot p^2 = 1\rangle ~\cong \\
\langle y, x_1,...,  x_m \mid y \cdot p^2 = 1\rangle$. 
\smallskip

 Upon applying the same trick $(k-1)$ times, we will have $y \cdot p = 1$ 
isomorphic to $y \cdot p^k = 1$.  $\Box$ \\

\noindent {\bf Acknowledgements}
\medskip

 We are grateful to A. Campbell and S. Kaliman for  helpful discussions. 
  The first author is grateful to the 
 Department of Mathematics of the 
University of Hong Kong for its warm hospitality during his visit
when most of  this work was done. \\

\noindent 
 Department of Mathematics, The City  College  of New York, New York, 
NY 10031 
\smallskip

\noindent {\it e-mail address\/}: ~shpil@groups.sci.ccny.cuny.edu 

\smallskip

\noindent {\it http://zebra.sci.ccny.cuny.edu/web/shpil} \\

\noindent Department of Mathematics, The University of Hong Kong, 
Pokfulam Road, Hong Kong 

\smallskip

\noindent {\it e-mail address\/}: ~yujt@hkusua.hku.hk 
\smallskip

\noindent {\it http://hkumath.hku.hk/\~~jtyu}

\end{document}